\documentclass[11pt]{article}

\usepackage{latexsym,mathrsfs}
\usepackage{amsmath,amssymb}
\usepackage{amsthm,enumerate,verbatim}
\usepackage{amsfonts}
\usepackage{graphicx}
\usepackage{algorithm}
\usepackage{algorithmic}
\usepackage{url}

\newtheorem{lemma}{Lemma}

\newtheorem{theorem}{Theorem}

\newtheorem{definition}{Definition}
\newtheorem*{definition*}{Definition}

\newtheorem{remark}{Remark}

\numberwithin{equation}{section}
\numberwithin{table}{section}
\numberwithin{figure}{section}

\DeclareMathOperator{\argmin}{argmin}

\newcommand {\mat}  [1] {\left[\begin{array}{#1}}
\newcommand {\rix}      {\end{array}\right]}

\title{On computing the distance to stability for matrices using linear dissipative Hamiltonian systems}
\author{Nicolas Gillis$^*$ \qquad Punit Sharma\thanks{The authors acknowledge the support of the ERC (starting grant n$^\text{o}$ 679515).
NG also acknowledges the support of the F.R.S.-FNRS (incentive grant for scientific research n$^\text{o}$ F.4501.16). } \\ % \thanks{xxx}
Department of Mathematics and Operational Research \\
Facult\'e Polytechnique, Universit\'e de Mons \\
Rue de Houdain 9, 7000 Mons, Belgium\\
 nicolas.gillis@umons.ac.be, Punit.SHARMA@umons.ac.be
}

\begin{document}

\maketitle

\begin{abstract}
In this paper, we consider the problem of computing the nearest stable matrix to an unstable one.
We propose new algorithms to solve this problem based on a reformulation using linear dissipative Hamiltonian systems:
we show that a matrix $A$ is stable if and only if it can be written as $A = (J-R)Q$, where $J=-J^T$, $R \succeq 0$
and $Q \succ 0$ (that is, $R$ is positive semidefinite and $Q$ is positive definite).
This reformulation results in an equivalent optimization problem with a simple convex feasible set.
We propose three strategies to solve the problem in variables $(J,R,Q)$:
(i) a block coordinate descent method, (ii) a projected gradient descent method, and
(iii) a fast gradient  method inspired from smooth convex optimization.
These methods require $\mathcal{O}(n^3)$ operations per iteration, where $n$ is the size of $A$.
We show the effectiveness of the fast gradient  method compared to the other approaches and to several state-of-the-art algorithms.
\end{abstract}


\section{Introduction}
In this paper, we focus on the continuous linear time invariant systems of the form
\[
\dot{x}(t)=Ax(t)+Bu(t),
\]
where $A\in \mathbb R^{n,n}$, $B\in \mathbb R^{n,m}$, $x$ is the state vector and $u$ is the input vector. Such
a system is stable if all eigenvalues of $A$ are in the \emph{closed} left half of the complex plane and all
eigenvalues on the imaginary axis are semisimple. Therefore the stability solely depends on
$A$, and the matrix $B$ that weights the inputs can be ignored to study stability. \\
For a given unstable matrix $A$, the problem of finding the smallest perturbation that stabilizes $A$, or,
equivalently finding the nearest stable matrix $X$ to $A$ is an important problem~\cite{ONV13}. More precisely,
we consider the following problem. For a given unstable matrix $A$, compute
\begin{equation}\label{eq:prob_def}
\inf_{X\in \mathbb S^{n,n}} {\|A-X\|}_F^2,
\end{equation}
where ${\|\cdot\|}_F$ denotes the Frobenius norm of a matrix and $\mathbb S^{n,n}$ is the set of
all stable matrices of size $n\times n$. This problem occurs for example in system identification where
one needs to identify a stable system from observations \cite{ONV13}.

The converse of problem~\eqref{eq:prob_def} is the stability radius problem, where a stable matrix $A$ is given and
one looks for the smallest perturbation that moves an eigenvalue outside the stability region~\cite{Bye88,HinP86}.
Both problems are nontrivial because
even a small perturbation on the coefficients of
the matrix may move the eigenvalues in any direction and the perturbed matrix may well
have eigenvalues that are far from those of $A$~\cite{ONV13}.
However, the nearest stable matrix problem appears to be more difficult since it requires to push all eigenvalues from the unstability region into the stability region while the stability radius problem only requires
 to move a single eigenvalue on the boundary of the stability region. \\
The various distance problems for matrices have been a topic of research in the numerical linear algebra community,
for example, matrix nearness problems~\cite{Hig88a}, the structured singular value
problem~\cite{PacD93}, the robust stability problem~\cite{Zho11}, the distance to bounded realness for Hamiltonian matrices~\cite{AlaBKMM11}, and the nearest defective matrix~\cite{Wil84}. \\
Another related problem is to find the closest stable polynomial to a given unstable one. This was addressed by Moses and Liu~\cite{MosL91}, where an algorithm using the alternating projection approach in Schur parameter space was
developed. But the technique developed in~\cite{MosL91} is limited and cannot be applied to other
types of systems. In~\cite{BurHLO06a}, authors stabilize fixed order controllers using nonsmooth, nonconvex optimization.
A MATLAB toolbox called HIFOO ($H_\infty$ fixed order optimization) was designed to solve fixed order stabilization
and local optimization problems~\cite{BurHLO06b}. In~\cite{DhaPV06}, authors stabilize transfer functions using a two step iterative procedure that guaranteed stable transfer function models from noisy data.

We note that in the literature, a stable matrix is sometimes considered to satisfy $\text{Re}(\lambda) < 0$ for all its
eigenvalues $\lambda$; see, e.g.,~\cite{Bye88,HinP86,ONV13}, .
 To avoid the confusion, we call such matrices asymptotically stable. The
set of all asymptotically stable matrices is open. This follows due to the fact that the eigenvalues of a matrix
depend continuously on its entries \cite{Ost60}. However, the set $\mathbb S^{n,n}$  is neither
open nor closed, because $A_\epsilon \notin \mathbb S^{n,n}$ for $\epsilon >0$, but $A \in \mathbb S^{n,n}$, where
\[
\underbrace{\mat{ccccc}
0 & 1 & 0 &\hdots &0\\
-1 & \epsilon & 0 &\hdots &0\\
0 & 0 & -1 &\hdots &0\\
\hdots&\hdots&\hdots&\hdots&\hdots\\
0 & 0 & 0 &\hdots &-1\\
\rix}_{=:A_\epsilon} \rightarrow
\underbrace{\mat{ccccc}
0 & 1 & 0 &\hdots &0\\
-1 & 0 & 0 &\hdots &0\\
0 & 0 & -1 &\hdots &0\\
\hdots&\hdots&\hdots&\hdots&\hdots\\
0 & 0 & 0 &\hdots &-1\\
\rix}_{ =:A},
\]
and
$B_\delta \in \mathbb S^{n,n}$ for $\delta < 0$, but $B \notin \mathbb S^{n,n}$, where
\[
\underbrace{\mat{ccccc}
0 & 1 & 0 &\hdots &0\\
0 & \delta & 0 &\hdots &0\\
0 & 0 & -1 &\hdots &0\\
\hdots&\hdots&\hdots&\hdots&\hdots\\
0 & 0 & 0 &\hdots &-1\\
\rix}_{=:B_\delta} \rightarrow
\underbrace{\mat{ccccc}
0 & 1 & 0 &\hdots &0\\
0 & 0 & 0 &\hdots &0\\
0 & 0 & -1 &\hdots &0\\
\hdots&\hdots&\hdots&\hdots&\hdots\\
0 & 0 & 0 &\hdots &-1\\
\rix}_{ =:B}.
\]
Further, the set $\mathbb S^{n,n}$ of stable matrices in~\eqref{eq:prob_def}
is highly nonconvex~\cite{ONV13} and therefore it is in general difficult to compute a global optimal
solution to problem~\eqref{eq:prob_def}.
% in a reasonable complexity.

Our work is mainly motivated by the work in~\cite{ONV13}, where a nearby stable
approximation $X$ of a given unstable system $A$ is constructed by means of successive convex approximations of the set of stable systems.
Our principle strategy for computing a nearby stable approximation to a given unstable matrix is to reformulate the highly nonconvex optimization problem~\eqref{eq:prob_def}
into an equivalent (non-convex) optimization problem with a convex feasible region onto which points can be projected relatively easily.
We aim to provide in many cases better approximations than the one obtained with the code from~\cite{ONV13} by
using the concept of linear dissipative Hamiltonian systems.

\textbf{Notation:} In the following, we denote $A\succ 0$ and $A\succeq 0$ if $A$ is symmetric positive definite
or symmetric positive semidefinite, respectively. The set $\Lambda(A)$ denotes the set of all eigenvalues of $A$.

\subsection{Dissipative Hamiltonian systems}
A dissipative Hamiltonian (DH) system in the linear time invariant case can be expressed as
\[
\dot{x}=(J-R)Qx,
\]
where the function $x\rightarrow x^TQx$ with $Q=Q^T \in \mathbb R^{n,n}$ positive definite describes
the energy of the system, $J=-J^T \in \mathbb R^{n,n}$ is the structure matrix that describes flux
among energy storage elements, and $R \in \mathbb R^{n,n}$ with $R=R^T \succeq 0$
is the dissipation matrix that describes energy dissipation in the system. DH systems
are special cases of port-Hamiltonian systems, which recently have received a lot attention
in energy based modeling; see, e.g.,~\cite{GolSBM03,Sch06,SchM13}. An important property of
DH systems is that they are stable, i.e., all eigenvalues of matrix $A=(J-R)Q$ are in
the closed left half of the complex plane and all eigenvalues on the imaginary axis are semisimple.
This follows from the fact that $Q$ is symmetric positive definite. Indeed, for any
nonzero vector $z$ one has
\begin{eqnarray*}
\mathrm{Re}{\left(z^*\big(Q^{\frac{1}{2}}AQ^{-\frac{1}{2}}\big)z\right)}&=&
\mathrm{Re}{\left(z^*\big(Q^{\frac{1}{2}}JQ^{\frac{1}{2}}-Q^{\frac{1}{2}}RQ^{\frac{1}{2}}\big)z\right )}\\
&=&-z^*Q^{\frac{1}{2}}RQ^{\frac{1}{2}}z \leq 0,
\end{eqnarray*}
since $R$ is positive semidefinite, where $*$ stands for the complex conjugate transpose
of a matrix or a vector. The semisimplicity of the purely imaginary eigenvalues
of $(J-R)Q$ follows from~{\rm \cite[Lemma 3.1]{MehMS16a}}.
 The various structured distances of a DH system from the
region of asymptotic stability have recently been studied in~\cite{MehMS16a} for the complex case and in~\cite{MehMS16b} for the real case. \\

This paper is organized as follows.
In section~2, we reformulate the nearest stable matrix problem using the notion of DH matrices.
We also provide several theoretical results necessary to obtain our reformulation.
In section~3, three algorithms are proposed to solve the reformulation.
In section~4, we present numerical experiments that illustrate the performance of our algorithms and
compare the results with several state-of-the-art algorithms.

\section{DH framework for checking stability}

In this section, we present a new framework based on dissipative Hamiltonian systems to
attack the nearest stable matrix problem~\eqref{eq:prob_def}. Our main idea is to
 reformulate the nonconvex optimization problem~\eqref{eq:prob_def}
into an equivalent optimization problem with a relatively simple convex feasible set.
In order to do this, let us define a DH matrix.
\begin{definition}
A matrix $A \in \mathbb R^{n,n}$ is said to be a {\em DH} matrix if $A=(J-R)Q$ for some $J,R,Q \in \mathbb R^{n,n}$ such that
$J^T=-J$, $R\succeq 0$ and $Q \succ 0$.
\end{definition}
Clearly from the previous section every DH matrix is stable.
In our terminology, {\rm \cite[Corollary 2]{BeaMX15}} implies that every stable matrix $A$ is similar to
a DH matrix, i.e., there exists $T$ nonsingular such that $T^{-1}AT=(J-R)Q$ for some
$J^T=-J$, $R\succeq 0$ and $Q \succ 0$. In fact we prove something stronger: a stable matrix itself is a
DH matrix, as shown in the following lemma.

\begin{lemma}\label{lem:stabisDH}
Every stable matrix is a {\em DH} matrix.
\end{lemma}
%\proof
\par\noindent{\bf Proof}. \ignorespaces
 Let $A$ be stable. By Lyapunov's theorem~\cite{LanT85}, there exists² $P \succ 0$ such that
 \begin{equation}\label{eq:lyapunov_cond}
 AP+P A^T \preceq 0.
 \end{equation}
Let us define
 \begin{equation} \label{eq:lyapunov_form_JRQ}
 J:=\frac{AP - (AP)^T}{2}, R:=-\frac{AP+(AP)^T}{2}, Q:=P^{-1}.
 \end{equation}
By~\eqref{eq:lyapunov_cond}
we get $R \succeq 0$ and $A=(J-R)Q$. This implies that $A$ is a DH matrix.
\space{\ \vbox{\hrule\hbox{\vrule height1.3ex\hskip0.8ex\vrule}\hrule}}\par
%\eproof

In~{\rm \cite[Lemma 3.1]{MehMS16a}}, authors provide several properties of DH matrices for $Q \succ 0$.
We now prove a similar result but for any $Q \succeq 0$, that is, allowing $Q$ to be singular.
\begin{lemma}\label{lem:stabJRQ}
Let $J,R,Q \in \mathbb R^{n, n}$ be such that $J^T=-J$, $R^T=R \succeq 0$, and $Q^T=Q \succeq 0$ singular.
Then
\begin{enumerate}
\item [(i)] All eigenvalues of $(J-R)Q$ are in the closed left half of the complex plane. Further,
$(\lambda,x) \in i\mathbb R\times \mathbb C^n\setminus \{0\}$ is an eigenpair of $(J-R)Q$
if and only if $RQx=0$ and $(JQ-\lambda I_n)x=0$,
where $i\mathbb R$ denotes the set of all purely imaginary complex numbers and
$I_n$ is the identity matrix of size $n$.
%$(J-R)Q$ has an
%eigenvalue $\lambda$ on the imaginary axis if and only if $RQx=0$ for some eigenvector $x$ of $JQ$
%with respect to eigenvalue $\lambda$.
\item [(ii)] All nonzero purely imaginary eigenvalues of the matrix $(J-R)Q$ are semisimple.
\end{enumerate}
\end{lemma}
%\proof
\par\noindent{\bf Proof}. \ignorespaces
 Proof of {\it{(i)}} is elementary and follows from the fact that $J^T=-J$, $R^T=R\succeq 0$, and
 $Q^T=Q \succeq 0$.
We prove {\it{(ii)}} by contradiction. Let $\lambda(\neq 0) \in i\mathbb R$ be an eigenvalue of
$(J-R)Q$. Suppose that $\lambda$ is not semisimple and that the set $\{x_0,x_1,\ldots,x_{m-1}\}$ forms a
Jordan chain of length $m$ at $\lambda$~\cite{HorJ85}, that is, $x_0 \neq 0$, and
\begin{eqnarray}\label{eq:lem1proof_1}
((J-R)Q-\lambda I_n)x_0&=&0,\nonumber\\
 ((J-R)Q-\lambda I_n)x_1&=&x_0,\\
((J-R)Q-\lambda I_n)x_2&=&x_1,\nonumber\\
\vdots\nonumber\\
((J-R)Q-\lambda I_n)x_{m-1}&=&x_{m-2}.\nonumber
\end{eqnarray}
Since $x_0$ is an
eigenvector of $(J-R)Q$ corresponding to the eigenvalue $\lambda$, from {\it{(i)}} we have that
\begin{equation}\label{eq:lem1proof_2}
(JQ-\lambda I_n)x_0=0\quad\text{and}\quad RQx_0=0.
\end{equation}
From~\eqref{eq:lem1proof_1}, $x_0$ and $x_1$ satisfy
\begin{equation}\label{eq:lem1proof_3}
((J-R)Q-\lambda I_n)x_1=x_0.
\end{equation}
Note that $\lambda \neq 0$ implies that $Qx_0 \neq 0$. Multiplying~\eqref{eq:lem1proof_3}
by $x^*Q$ from the left hand side we obtain that
\[
x_0^*Q((J-R)Q-\lambda I_n)x_1=x_0^*Qx_0.
\]
This implies that
\begin{equation}\label{eq:lem1proof_4}
-x_1^*Q(JQ-\lambda I_n)x_0 -x_1^*QRQx_0=x_0^*Qx_0,
\end{equation}
because $Q \succeq 0$, $R\succeq 0$, and $J^T=-J$. Using~\eqref{eq:lem1proof_2}
in~\eqref{eq:lem1proof_4}, we get that
$ x_0^*Qx_0=0$. This implies that $x_0=0$ as $Q \succeq 0$ and $Qx_0 \neq 0$
which is a contradiction since the generator of the Jordan chain
$x_0$ is never zero.
Therefore $\lambda$ is a semisimple eigenvalue of $(J-R)Q$.
\space{\ \vbox{\hrule\hbox{\vrule height1.3ex\hskip0.8ex\vrule}\hrule}}\par
%\eproof

\begin{remark}{\rm
We note that if $Q$ is positive definite, then we have $JQ=Q^{-\frac{1}{2}}(Q^{\frac{1}{2}}JQ^{\frac{1}{2}})Q^{\frac{1}{2}}$, i.e.,
$JQ$ is similar to a real skew-symmetric matrix $Q^{\frac{1}{2}}JQ^{\frac{1}{2}}$. Therefore, all
eigenvalues of $JQ$ are purely imaginary and semisimple, and as a consequence
of~{\rm \cite[Lemma 3.1]{MehMS16a}} all purely imaginary eigenvalues of $(J-R)Q$ are semisimple.
If $Q \succeq 0$ and singular, then by
Lemma~\ref{lem:stabJRQ} all nonzero purely imaginary eigenvalues of $(J-R)Q$ are semisimple but
zero eigenvalues of $(J-R)Q$ can still be non semisimple and hence the matrix is not stable. However,
this matrix is on the boundary of the set of unstable matrices. This gives us a motivation to look for the
nearest matrix $(J-R)Q$, with $J^T=-J,R\succeq 0$, and $Q \succeq 0$ that is either stable or on the boundary of the
set of unstable matrices.
}
\end{remark}

In view of Lemma~\ref{lem:stabisDH} and  Lemma~\ref{lem:stabJRQ}, a matrix is stable if and only if
it is a DH matrix. Thus the set $\mathbb S^{n,n}$ of stable matrices can be characterized as the set of
DH matrices, i.e.,
\[
\mathbb S^{n,n}=\left\{(J-R)Q\in \mathbb R^{n,n}\,|\,J^T=-J,\,R\succeq 0, \,Q \succ 0\right \}.
\]
This yields the following equivalent reformulation of the nonconvex problem~\eqref{eq:prob_def}:
\begin{equation} \label{eq:reformulate_prob}
\inf_{X\in \mathbb S^{n,n}}{\|A-X\|}_F^2=\inf_{J=-J^T,\, R \succeq 0,\, Q \succ 0}{\left\|A - (J-R)Q\right\|}_F^2.
\end{equation}
DH characterization of the set $\mathbb S^{n,n}$ also gives an alternative way to see that the set $\mathbb S^{n,n}$
is neither open (because of the constraint $R\succeq 0$) nor closed (because of the constraint $Q \succ 0$).
%
%is an open set, i.e., marginally stable matrices may not belong to the set. This happens when $Q$ is positive
%semidefinite and singular, and zero is not a semisimple eigenvalue of $(J-R)Q$, hence the search for an
%infimum in~\eqref{eq:prob_def}, or equivalently in~\eqref{eq:reformulate_prob}.
%
Further, we have
\begin{equation}\label{eq:reformulate_process_1}
\inf_{J=-J^T,\, R \succeq 0,\, Q \succ 0}f(J,R,Q)=\inf_{J=-J^T,\, R \succeq 0,\, Q \succeq 0}f(J,R,Q),
\end{equation}
where $f(J,R,Q)={\left\|A - (J-R)Q\right\|}_F^2$. Note that the set 
\[
\left\{(J,R,Q)\,|\,J,R,Q \in \mathbb R^{n,n},\,J^T=-J,\,R\succeq 0, \,Q \succeq 0\right \} 
\] 
is closed but not bounded. Therefore the infimum in the right hand side of~\eqref{eq:reformulate_process_1}
may not be attained.
As a summary we have the following desired reformulation of our problem~\eqref{eq:prob_def}.
\begin{theorem}\label{thm:main_reformulation}
Let $A \in \mathbb R^{n,n}$. Then,
\begin{equation} \label{reformulation}
\inf_{X\in \mathbb S^{n,n}}{\|A-X\|}_F^2=\inf_{J=-J^T,\, R \succeq 0,\, Q \succeq 0}{\left\|A - (J-R)Q\right\|}_F^2.
\end{equation}
\end{theorem}

In the next section, we will attack the nearest stable matrix problem~\eqref{eq:prob_def} using the reformulation from
Theorem~\ref{thm:main_reformulation}, i.e., by trying to solve
\begin{equation}\label{nearportham}
\inf_{J=-J^T,\, R \succeq 0,\, Q \succeq 0}{\left\|A - (J-R)Q\right\|}_F^2.
\end{equation}
An advantage of this reformulated optimization problem over~\eqref{eq:prob_def} is that the
feasible domain of~\eqref{nearportham} is convex and that it is relatively easy to project onto it (see Section~\ref{sec:pgrad}).
Hence it will be easy to use standard optimization schemes to solve it.

If the matrix $A$ is stable, then by Lemma~\ref{lem:stabisDH} it can be written as $A = (J-R)Q$ with
$J=-J^T$, $R \succeq 0$ and $Q \succ 0$, hence $AQ^{-1} = J-R$. In that case, we can solve the following
system to recover $(J,R,Q)$: denoting $P = Q^{-1}$,
\[
AP = J-R, \; P \succ 0, \; R \succeq 0, \; J = -J^T.
\]
This is interesting because it provides a new (convex) way to check whether a matrix is stable, checking whether the above system is feasible. Moreover, if $A$ is not stable so that the above system is infeasible, we can solve instead
\begin{equation} \label{nearporthaminv}
\inf_{J=-J^T,\, R \succeq 0,\,P \succeq I_n} {\| AP - (J-R) \|}_F^2,
\end{equation}
which provides an approximate solution to~\eqref{nearportham} using $(J-R)P^{-1}$ as a stable approximation of $A$.
This solution could be used as an initialization for nearest stable matrix algorithms.
From the standard stability formulation~\eqref{eq:lyapunov_cond}, as far as we know, it is not possible to extract a stable approximation from an unstable matrix. This is another advantage of our formulation.

\begin{remark}[Uniqueness]{\rm \label{remunique}
The decomposition of a stable matrix $A$ as $A = (J-R)Q$, where $J^T=-J, R\succeq 0$, and $Q \succ 0$, is non-unique.
First, there is always a scaling degree of freedom since $A = (\alpha J- \alpha R) (Q/\alpha)$ for any scalar $\alpha > 0$.
In the numerical algorithms, we will impose that ${\|J-R\|}_2={\|Q\|}_2$ to remove this degree of freedom;
see the discussion in Section~\ref{linesearch}.
Second, in view of~\eqref{eq:lyapunov_cond} and~\eqref{eq:lyapunov_form_JRQ}, the non-uniqueness of $(J-R)Q$ can be (partly) characterized with the matrix $P \succ 0$ that certifies the stability of $A$, which is not unique~\cite{BoyGFB94}.
Intuitively, the set of stable matrices is an $n^2$-dimensional space while the representation
$(J-R)Q$ has $\frac{3}{2} n^2 + \frac{n}{2}$ degrees of freedom (since $R$ and $Q$ are symmetric, and $J$ is skew-symmetric).
Replacing $J$, $R$ and $Q$ in terms of the variable $P \succ 0$
 according to~\eqref{eq:lyapunov_form_JRQ} leads to the formulation
\[
\min_{X, P \succ 0} || A - X ||_F \; \text{ such that } \; X P + P X^T \preceq 0,
\]
which is difficult to solve numerically~\cite{ONV13} (highly non-linear constraints).
In fact, a key contribution of this paper is the reformulation~\eqref{eq:reformulate_process_1} of the nearest stable matrix problem in term of DH matrices that will allow us to derive numerically efficient algorithms;
see Section~\ref{algoDHmat}.
%which is difficult to solve numerically (highly non-linear constraints).
%However, this would not be a desirable formulation for our numerical algorithms since we would need to deal with a highly non-linear problem, in particular involving $P^{-1}$.
Characterizing precisely the non-uniqueness of this decomposition (and possibly taking advantage of it in a numerical algorithm)
%seems non-trivial and out of
 is a direction for further research.
}
\end{remark}

\section{Algorithmic solutions to the nearest stable matrix problem} \label{algoDHmat}
In this section, we propose three algorithms to compute approximate solutions to our reformluation of the nearest stable matrix problem~\eqref{nearportham} in variables $J,R$ and $Q$.

\subsection{Block coordinate descent method} \label{sec:bcd}

Popular and simple schemes in optimization are the block-coordinate descent (BCD) methods: at each step,
 fix a subset of the variables and optimize over the remaining variables; see, e.g.,~\cite{W2015,STY16} for recent surveys.
For our problem~\eqref{nearportham}, there is a natural partition of the variables in order to make the subproblems convex:
$(J,R)$ vs.\@ $Q$.
Hence, we propose the following simple scheme:
\begin{itemize}
\item Initialize $(J,R,Q)$ such that $J=-J^T$, $R \succeq 0$ and $Q \succeq 0$.
\item Iteratively update $(J,R)$ and $Q$:
\begin{itemize}
\item []{\em subproblem-1}: fix $Q$ and optimize over $(J,R)$, that is, solve
\[
\min_{J=-J^T,\,R \succeq 0} {\|A - (J-R)Q\|}_F^2.
\]
\item []{\em subproblem-2}: fix $(J,R)$ and optimize over $Q$, that is, solve
\[
\min_{Q \succeq 0} {\|A - (J-R)Q\|}_F^2.
\]
\end{itemize}
\end{itemize}
Using dedicated semidefinite programming (SDP) solvers such as SDPT3~\cite{toh1999sdpt3} or SeDuMi~\cite{sturm1999},
this can be very easily implemented; especially if combined with a modeling system for convex optimization like CVX~\cite{grant2008cvx}.
However, (most) SDP solvers use interior methods which are computationally expensive when $n$ increases
(as they perform Newton steps at their core). Moreover, we have to solve these problems many times, and it is not necessary to solve them up to high precision.
We have therefore implemented fast gradient  methods to solve these subproblems: these are optimal first-order methods for smooth convex optimization; see~\cite[p.90]{nes04} and Appendix~\ref{appA}. 

This method works relatively well, especially to decrease the objective function initially. However, when getting closer to a stationary point, it has a zig-zaging behavior and converges rather slowly.
We note that the method is guaranteed to converge to a stationary point of~\eqref{nearportham}
since it is a two-block coordinate descent method~\cite{GS00}.

Since we use a first-order method to solve the subproblems, the total computational cost of this approach is $O(K n^3)$ operations per iteration, where $K$ is the number of inner iterations. In our implementation, we used $K=1000$ along with a stopping criterion based on the evolution of the iterates:
if the distance between two iterates is smaller than 1\% of the distance between the first two iterates,
we stop the inner iterations --this means that the improvement is negligible compared to the first iteration hence it is profitable to start optimizing over the other set of variables.

\subsection{Projected Gradient Descent} \label{sec:pgrad}
One of the most well-known and straightforward algorithm one can implement for~\eqref{nearportham} is projected gradient
descent. For a given matrix $Z\in \mathbb R^{n,n}$, we recall that the projection
$\mathcal P_{\bar S}(Z)$ of $Z$ onto the set of skew-symmetric matrices $\bar S$ is
given by
\begin{equation}\label{eq:sec_3_eq_1}
\mathcal P_{\bar S}(Z)=\frac{Z-Z^T}{2},
\end{equation}
and we have
\begin{equation}\label{eq:sec_3_eq_2}
\min_{J^T=-J}{\|Z-J\|}_F={\|Z-\mathcal P_{\bar S}(Z)\|}_F={\left\|\frac{Z+Z^T}{2}\right\|}_F.
\end{equation}
Similarly, the projection
$\mathcal P_{ \succeq }(Z)$ of $Z$ onto the cone of semidefinite matrices is given by
\begin{equation}\label{eq:sec_3_eq_3}
\mathcal P_{ \succeq }(Z)=U\left(\max{(\Gamma,0)}\right)U^T,
\end{equation}
where $U \Gamma U^T$ is an eigenvalue decomposition of the symmetric matrix $\frac{Z+Z^T}{2}$,
and we have
\begin{eqnarray}\label{eq:sec_3_eq_4}
\min_{R\succeq 0}{\|Z-R\|}_F^2&=&{\|Z-\mathcal P_{\succeq }(Z)\|}_F^2\\
&=&{\|\mathcal P_{\bar S}(Z)\|}_F^2+\sum_{\lambda \in \Lambda(\Gamma),\lambda < 0}\lambda^2,\nonumber
\end{eqnarray}
see, e.g.,~\cite{Hig88b}. This follows directly from univariance of ${\|\cdot\|}_F$ to orthogonal transformations.

 In order to simplify the description of our algorithms,
we define
\[
\mathcal{D} := \left\{ D \in \mathbb{R}^{n,n} \, |~ D = J - R,~ J = -J^T,~ R \succeq 0 \right\}.
\]
It is interesting to note that it is as easy to project a given matrix $Z$ onto $\mathcal{D}$,
as it is equivalent to project onto the set of skew-symmetric matrices, $\bar{\mathcal{S}}$, and
 the set of positive semidefinite matrices, $\mathcal{S}_{\succeq }$, separately.
More precisely, we have the following lemma.

\begin{lemma}\label{lem:sec_2_lem}
Let $Z \in \mathbb R^{n,n}$, then
\[
\min_{D\in \mathcal D}{\|Z-D\|}_F^2={\|Z-(\mathcal P_{\bar S}(Z)-\mathcal P_{\succeq}(-Z))\|}_F^2,
%=\sum_{\lambda \in \Lambda(\Gamma),\lambda > 0}\lambda^2,
\]
where $\mathcal P_{\bar S}(Z)$ and $\mathcal P_{\succeq }(-Z)$ are as defined in~\eqref{eq:sec_3_eq_1}
and~\eqref{eq:sec_3_eq_3}, respectively.
\end{lemma}
%\proof
\par\noindent{\bf Proof}. \ignorespaces
Observe that
\begin{eqnarray}\label{eq:skewsym_def}
&\min_{D\in \mathcal D}{\|Z-D\|}_F^2=\min_{J^T=-J,\,R\succeq 0}{\|Z-(J-R)\|}_F^2 \nonumber\\
&= \min_{R\succeq 0} \left(\min_{J^T=-J}{\|(Z+R)-J\|}_F^2\right) \nonumber\\
&= \min_{R\succeq 0}{\left\|(Z+R)-\mathcal P_{\bar S}(Z)\right\|}_F^2
= \min_{R\succeq 0}{\left\|\frac{Z+Z^T}{2}+R\right\|}_F^2,\nonumber\\
\end{eqnarray}
where the third equality is due to~\eqref{eq:sec_3_eq_2} because
\[
\mathcal P_{\bar S}(Z+R)=\frac{(Z+R)-(Z+R)^T}{2}=\frac{Z-Z^T}{2}=\mathcal P_{\bar S}(Z).
\]
Thus the result follows by using~\eqref{eq:sec_3_eq_4} in~\eqref{eq:skewsym_def}.
\space{\ \vbox{\hrule\hbox{\vrule height1.3ex\hskip0.8ex\vrule}\hrule}}\par
%\eproof
%
Therefore, from Lemma~\ref{lem:sec_2_lem}, the projection of $Z$ onto $\mathcal D$ is given by
\[
\mathcal P_{\mathcal D}(Z)=\mathcal P_{\bar S}(Z)-\mathcal P_{\succeq }(-Z).
\]
In the following sections, by $D$ we mean a matrix of the form $J-R$ with $J^T=-J$ and $R\succeq 0$, and we denote
 $f(D,Q) := \frac{1}{2} {\| A - DQ\|}_F^2$, so our aim is to solve $\min_{D \in \mathcal{D},  Q \succeq 0} f(D,Q)$.
 We implemented a projected gradient decent algorithm for the
 nearest stable matrix problem~\eqref{nearportham}; see Algorithm~\ref{alg:pgrad}.
The computational cost of Algorithm~\ref{alg:pgrad} requires $O(n^3)$ operations per iteration.

\algsetup{indent=2em}
\begin{algorithm}[ht!]
\caption{Projected Gradient Descent (Grad) for the Nearest Stable Matrix Problem~\eqref{nearportham}} \label{alg:pgrad}
\begin{algorithmic}[1]
\REQUIRE $A \in \mathbb{R}^{n,n}$, {\ttfamily maxiter}.
    \ENSURE $J \in \bar{\mathcal{S}}$, $R \succeq 0$ and $Q \succeq 0$ such that ${\|A-(J-R)Q\|}_F$ is minimized.
    %\STATE \COMMENT{Initialization}
    \STATE Choose some initial matrices $J \in \bar{\mathcal{S}}$, $R \succeq 0$ and $Q \succeq 0$ (see Section~\ref{init}), set $D = J-R$.
	  \FOR{$k = 1 :$ {\ttfamily maxiter}}
        \STATE Compute the search direction (opposite of the gradient):
				\[
				\Delta_D  = -\nabla_D f(D,Q) = (A - DQ)Q^T,
				\]
				\[
				\Delta_Q  = -\nabla_Q f(D,Q) = D^T(A - D Q).
				\]
		\STATE Compute a steplength $\delta$ (see Section~\ref{linesearch}) and the next iterate:
\[
D = \mathcal{P}_{\mathcal{D}} \left( D + \delta \Delta_D \right),
\]
\[
Q = \mathcal{P}_{\succeq } \left( Q + \delta \Delta_Q \right).
\]
\STATE Scale $D,Q$ so that $\lambda_{\max}(QQ^T) = \lambda_{\max}(D^TD)$; see Section~\ref{linesearch}.
    \ENDFOR
\STATE $J = \frac{D-D^T}{2}$, $R =  \frac{-D-D^T}{2}$.
\end{algorithmic}
\end{algorithm}

\subsubsection{Line Search and Scaling}  \label{linesearch}
As noted in the previous section, if $D$ is fixed (resp.\@ $Q$), the problem is convex in $Q$ (resp.\@ in $D$).
If we would apply a gradient method on these subproblems,
a possible choice for the steplength that would guarantee convergence is $1/L$ where $L$ is the Lipschitz constant of the gradient of the objective function~\cite{nes04}.
The Lipschitz constant of the gradient $\nabla_Q f(D,Q)$ of $f$ with respect to $Q$  (for fixed $D$) is given by $\lambda_{\max}( D^T D )$
while the Lipschitz constant of the gradient $\nabla_D f(D,Q)$ of $f$ with respect to $D$ (for fixed $Q$)  is given by $\lambda_{\max}( QQ^T )$,
where we used $\lambda_{\max}(Z)$ to denote the largest eigenvalue of a matrix $Z$.

Therefore, it makes sense to scale $D$ and $Q$ such that  $L = \lambda_{\max}( D^T D ) = \lambda_{\max}( QQ^T )$ while
choosing an initial steplength $\delta = 1/L$. Note that this allows to remove the scaling degree of freedom since this imposes
 ${\|J-R\|}_2={\|Q\|}_2$; see also Remark~\ref{remunique}.
%This would allow to maximize the smallest steplength between $D$ and $Q$.
In order to avoid to having to compute the maximum eigenvalues of $D^T D$ and $QQ^T$ from scratch at each step,
we use a few steps of the power method to update the initial value (since $Q$ and $D$ do not change too much between two iterations).
We combined this with a backtracking line search: if the objective function has not decreased, the step is divided by a fixed constant larger than one until decrease is achieved. It turns out that in most cases, especially when getting closer to stationary points,
the steplength of $1/L$ allows to decrease the objective function.

\subsection{Fast Projected Gradient Descent} \label{sec:fgm}
For smooth convex optimization problems, the fast gradient  method (FGM), also known as Nesterov's method or the accelerated gradient method,
guarantees the objective function to decrease to the optimal value at a rate of $\mathcal{O}\left( \frac{1}{k^2} \right)$ where $k$ is the iteration number (or linearly with factor $(1-\sqrt{\mu/L})$ for strongly convex function with parameter $\mu$ and whose gradient has Lipschitz constant $L$) as opposed to gradient method converging at a rate of $O(1/k)$ (or linearly with factor $(1-{\mu/L})$ in the strongly convex case).
This method is an optimal first-order method for smooth convex optimization,
that is, no method can achieve asymptotically faster convergence using the same (first-order) information~\cite{nes04}.
In a nutshell, the idea behind FGM is to use the information from previous iterates (as opposed to using the information only of the previous iterate in gradient descent) similarly as conjugate gradient methods;
we refer the reader to~\cite{nes04} for more details and discussions.

For non-convex problems, the behavior of FGM is more difficult to analyze and, as far as we know, is not very well understood.
Ghadimi and Lan~\cite{GL16} have shown recently that FGM can be used successfully in the non-convex and unconstrained case, and achieves the best known rate of convergence in that case also.
However, there is straightforward practical way to ensure convergence using a simple safety procedure: if the FGM does not decrease the objective function, it can be reinitialized by performing a standard gradient step (note that this idea is also sometimes used for FGM in the convex case as it has been observed to work sometimes better, the reason being that FGM
does not decrease the objective function monotonically).
Moreover, doing so every fixed number of iterations allows to apply standard convergence results of the projected gradient method in non-convex optimization.

\algsetup{indent=2em}
\begin{algorithm}[ht!]
\caption{Fast Projected Gradient Method (FGM) for the Nearest Stable Matrix Problem~\eqref{nearportham}} \label{alg:fgm}
\begin{algorithmic}[1]
\REQUIRE $A \in \mathbb{R}^{n,n}$, {\ttfamily maxiter}, $0 < \alpha_{1} < 1$.
    \ENSURE $J \in \bar{\mathcal{S}}$, $R \succeq 0$ and $Q \succeq 0$ such that ${\|A-(J-R)Q\|}_F$ is minimized.
    \STATE Choose some initial matrices $J \in \bar{\mathcal{S}}$, $R \succeq 0$ and $Q \succeq 0$ (see Section~\ref{init}).
		\STATE $D = J-R$.
		\STATE $Y = D$, $Z = Q$.
    \FOR{$k = 1 :$ {\ttfamily maxiter}}
		\STATE $D_p = D$, $Q_p = Q$.
		    \STATE Compute the search direction (opposite of the gradient):
				\[
				\Delta_Y  = -\nabla_Y f(Y,Q) = (A - YQ)Q^T,
				\]
				\[
				\Delta_Z  = -\nabla_Z f(D,Z) = D^T(A - D Z).
				\]
		\STATE Compute a steplength $\delta$ (see Section~\ref{linesearch}) and the next iterate:
\[
D = \mathcal{P}_{\mathcal{D}} \left( Y + \delta \Delta_Y \right),
\]
\[
Q = \mathcal{P}_{\succeq } \left( Z + \delta \Delta_Z \right).
\]
\IF{no steplength allows the decrease of the objective function}
\STATE Restart the fast gradient : $Y = D$, $Z = Q$, $\alpha_{k+1} = \alpha_1$.
\ELSE{}
\STATE Update $Y$ and $Z$:
\[
Y = D + \beta_k (D-D_p),~Z = Q + \beta_k (Q-Q_p),
\]
 where $\beta_k =  \frac{\alpha_{k} (1-\alpha_{k})}{\alpha_{k}^2 + \alpha_{k+1}}$
with
$\alpha_{k+1} \geq 0$ s.t.
$\alpha_{k+1}^2 = (1-\alpha_{k+1}) \alpha_{k}^2$.
\ENDIF
\STATE Scale $Y,Z,D,Q$ so that $\lambda_{\max}(QQ^T) = \lambda_{\max}(D^TD)$; see Section~\ref{linesearch}.
\ENDFOR
\STATE $J = \frac{D-D^T}{2}$, $R =  \frac{-D-D^T}{2}$.
\end{algorithmic}
\end{algorithm}

For the fast gradient  version, we use the same line search as for the gradient method. However, in this case, because the algorithm is not guaranteed to be monotonously decreasing the objective function, the next iterate might not decrease the objective function for any steplength. In that case, as mentioned above, we reinitialize the fast gradient  scheme which amounts to perform a gradient step at the next iteration.

A possible explanation of the good performances of this fast gradient method designed for nonconvex problem is that
(1)~the feasible set is convex, and
(2)~the objective function can be locally well approximated with a convex function:
\begin{align*}
& {\| A - (D+\Delta_D)(Q+\Delta_Q) \|}_F   \\
 & \qquad = {\| A - D Q - \Delta_D Q - D \Delta_Q - \Delta_D \Delta_Q \|}_F \\
 &  \qquad \approx {\| A - D Q - \Delta_D Q - D \Delta_Q \|}_F,
\end{align*}
for small $(\Delta_D, \Delta_Q)$ since ${\|\Delta_D \Delta_Q \|}_F$ will be small compared to the other terms.
In other words, although~\eqref{nearportham} is non-convex,
it is close to a convex optimization problem in the neighborhood of any point $(D,Q)$.
In fact, rather surprisingly, we observed the following
\begin{quote}
%Except sometimes during the first iterations:
 (1)~the FGM almost never needs to restart, and \\
 (2)~the steplength of $1/L$ is almost always acceptable (making the search for a steplength inexpensive in most cases).
\end{quote}
In our experiments, FGM needed to restart or the steplength 1/L was not acceptable mostly during the first iterations.

The computational cost per iteration of FGM is (almost) the same as the gradient descent method, with $O(n^3)$ operations.

\subsection{Initialization} \label{init}
A natural way to initialize our algorithm is to use $Q = I_n$, where $I_n$ is the identity matrix of size $n$.
In that case, by Lemma~\ref{lem:sec_2_lem}, the optimal $J$ and $R$ are given by
\[
J = \frac{A-A^T}{2} \quad \text { and } \quad R = \mathcal{P}_{\succeq }\left(\frac{-A-A^T}{2}\right) .
\]
It seems this initialization works remarkably well. In particular, it allows the algorithms to identify better local minima
than using the strategy based on the formulation~\eqref{nearporthaminv}, although we observed it has in most cases a higher initial error.

\section{Numerical Experiments}

In this section, we will compare our three approaches to solve~\eqref{nearportham} (namely, BCD from Section~\ref{sec:bcd}, Grad from Section~\ref{sec:pgrad}, and FGM from Section~\ref{sec:fgm}) with the following two algorithms:
\begin{itemize}
\item SuccConv. A method from~\cite{ONV13} that uses successive convex approximation,
kindly made available to us by Fran\c{c}ois-Xavier Orban de Xivry.
\item BFGS. A method based on the BFGS code HANSO v.2.1 available from

\url{http://www.cs.nyu.edu/overton/software/hanso/}

and kindly made available to us by Michael Overton~\cite{LO09}. Recall that BFGS is a quasi-Newton method with a particular rank-two correction of the approximation of the Hessian at each iteration.
Note that this algorithm is based on a penalty method hence does not guarantee all iterates to be stable.
For $n \geq 50$, we used the limited memory variant (using 50 rank-one factors to estimate the Hessian).
Moreover, the algorithms terminates when it cannot find a descent direction (note that we tried to reinitialize the algorithm from the last iterate but it would terminate again, after one iteration).
\end{itemize}

We will compare these approaches on four types of matrices:
\begin{itemize}
\item Type 1. The matrices proposed in~\cite{ONV13}:
the entries on the subdiagonal are equal to one and the entry at position $(1,n)$ is equal to $-0.1$.
The eigenvalues are located on a circle around the origin; see Figure~\ref{type1n10}.
The best known approximation on the boundary of the set of stable matrices is the same matrix where the entry at position $(1,n)$ is set to zero, with error 0.1.
\item Type 2. Grcar matrices of order $k$ are a banded Toeplitz matrix with its subdiagonal set to $-$1 and both its main and $k$ superdiagonals set to 1. These matrices have all their eigenvalues in the right complex plane. This type of matrices was kindly recommended by Nicola Guglielmi during his talk `Matrix stabilization using differential equations' at the 9th Workshop on Structural Dynamical Systems: Computational Aspects (SDS2016); see

\url{https://sites.google.com/site/workshopsds2016/programm}.

We will use the same value of $k$, namely $k=3$.
These matrices have all their eigenvalues in the right half of the complex plane.
\item Type 3. Randomly generated matrices where each entry follows a Gaussian distribution N(0,1). These matrices have in average half their eigenvalues in the left complex plane and half in the right.
\item Type 4. Randomly generated matrices where each entry follows a uniform distribution on the interval [0,1]. Except for an eigenvalue with large  positive part (Perron-Frobenius theorem), the others are evenly spread on the right and left halves of the complex plane.
\end{itemize}

For each type of matrices, we generate a matrix of size 10, 20, 50 and 100. We run the algorithm for
at most 20, 100, 300 and 600 seconds respectively and record the final error.
For all algorithms, we use the initialization described in Section~\ref{init} because it performed consistently better than random initializations and than the initialization obtained from solving~\eqref{nearporthaminv}.
For $n = 100$, we do not run SuccConv because one iteration takes more than 10 minuts (it takes about 2 minuts for $n=60$).

The Matlab code is available from \url{https://sites.google.com/site/nicolasgillis/}.  All tests are preformed using Matlab on a laptop Intel CORE i5-3210M CPU @2.5GHz 2.5GHz 6Go RAM.
Note that FGM can be applied to matrices of size up to a thousand on such a machine, although we do not provide numerical experiments for $n$ larger than 100.

Tables~\ref{results10}, \ref{results20}, \ref{results50} and \ref{results100}
%give the final error ${\|A-X\|}_F$ obtained by each algorithm for each input matrix. Each table
display the final error ${\|A-X\|}_F$ obtained by the matrix $X$ generated by each algorithm and, in brackets, the number of iterations performed. Bold indicates the error of the best solution found.

\begin{figure}[h!] 
\begin{center}
\includegraphics[width=0.9\textwidth]{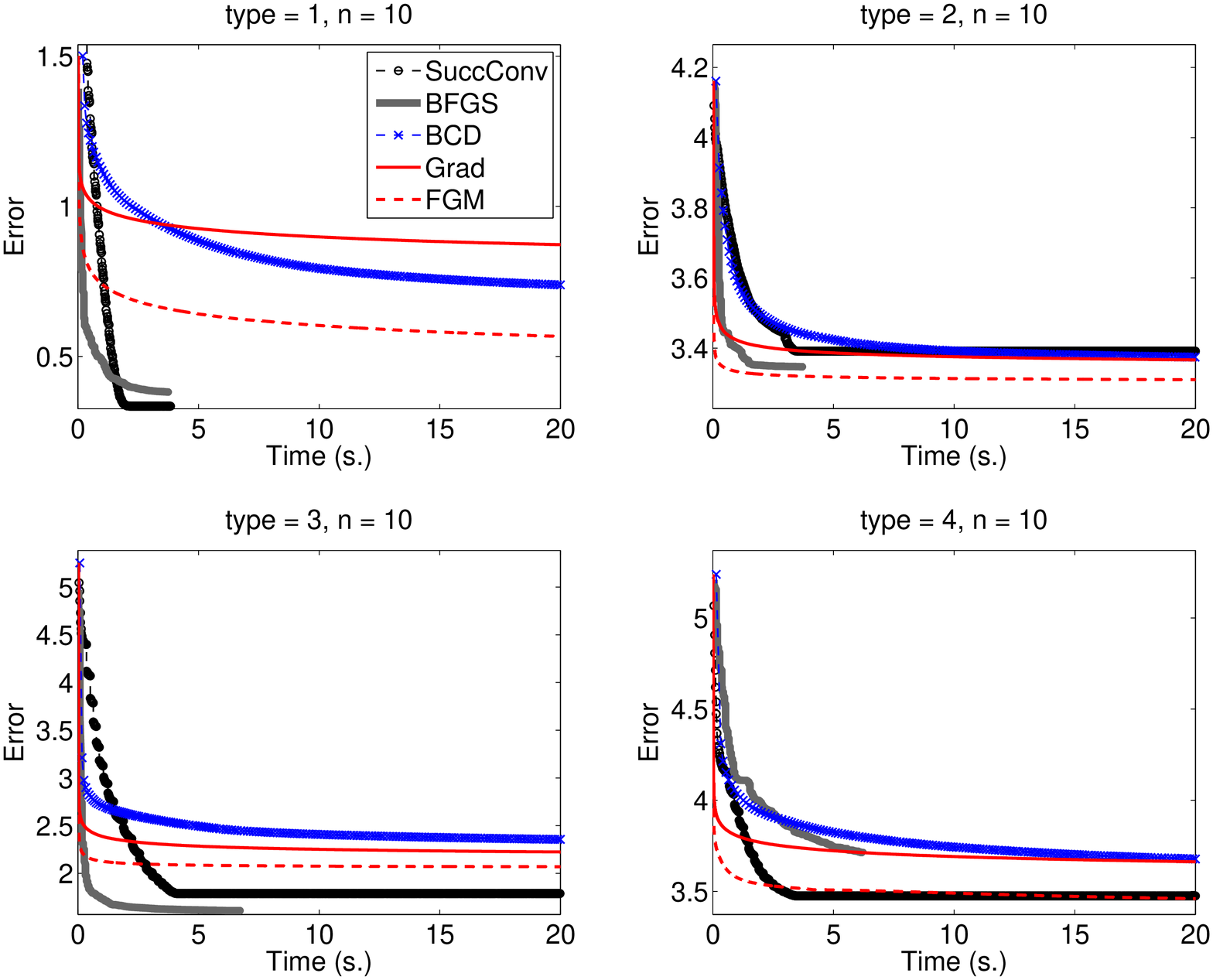}
\caption{Evolution of the error $\|A-X\|_F$ for the different algorithms for matrices $A$ of size 10. \label{matn10}}
\end{center}
\end{figure}

\begin{figure}[h!]
\begin{center}
\includegraphics[width=0.9\textwidth]{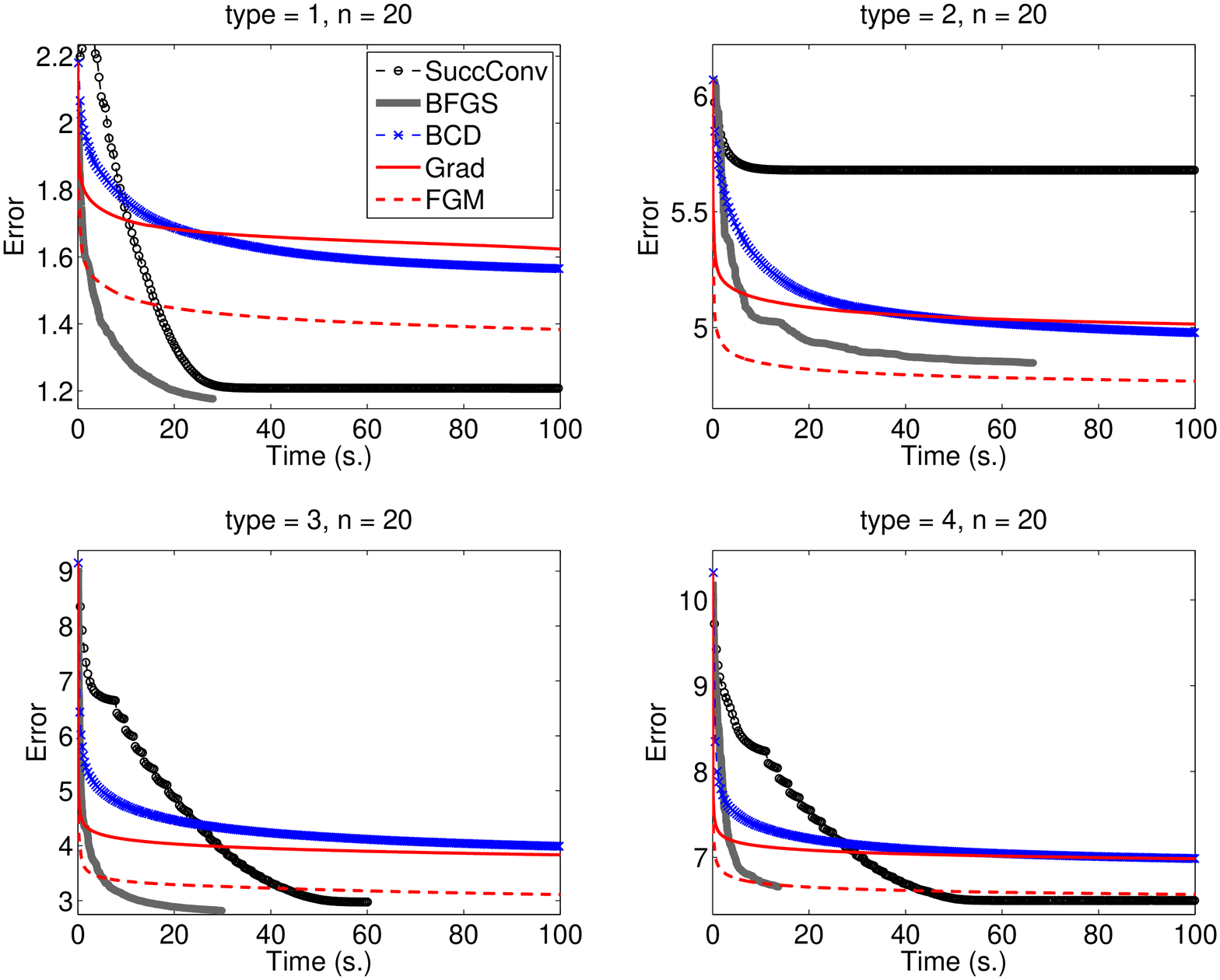}
\caption{Evolution of the error $\|A-X\|_F$ for the different algorithms for matrices $A$ of size 20. \label{matn20}}
\end{center}
\end{figure}

\begin{figure}[h!]
\begin{center}
\includegraphics[width=0.9\textwidth]{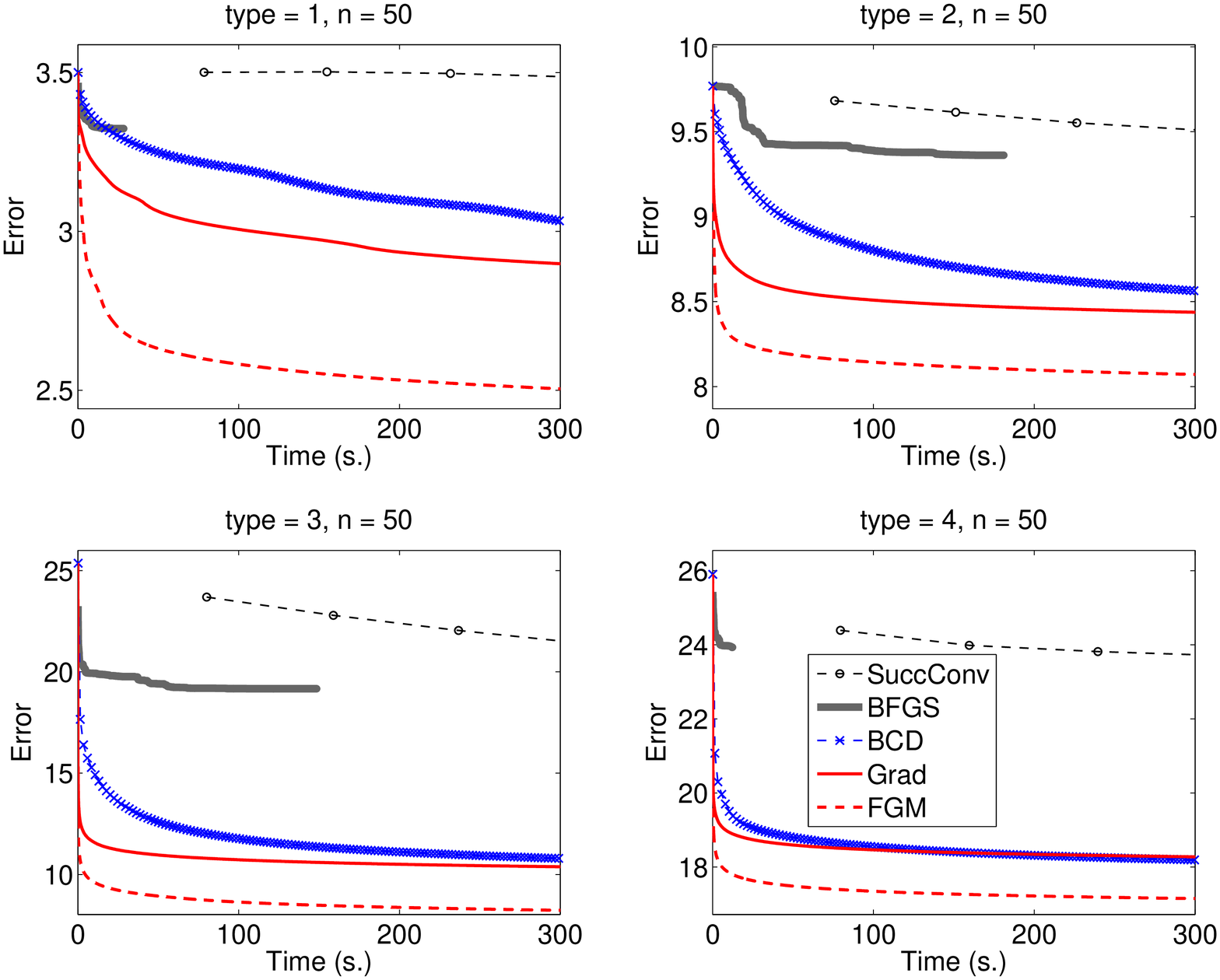}
\caption{Evolution of the error $\|A-X\|_F$ for the different algorithms for matrices $A$ of size 50. \label{matn50}}
\end{center}
\end{figure}

\begin{figure}[ht!]
\begin{center}
\includegraphics[width=0.9\textwidth]{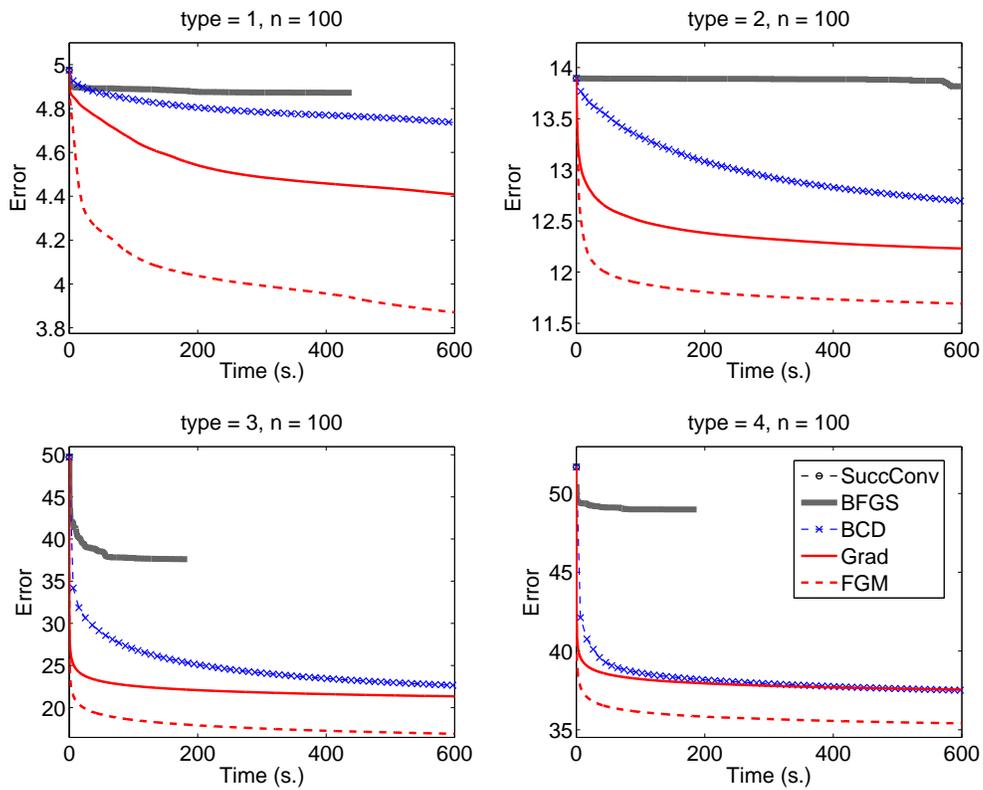}
\caption{Evolution of the error $\|A-X\|_F$ for the different algorithms for matrices $A$ of size 100. \label{matn100}}
\end{center}
\end{figure}

Figures~\ref{matn10}, \ref{matn20}, \ref{matn50} and \ref{matn100} display the evolution of the error ${\|A-X\|}_F$ for each algorithm for each input matrix.
 
\newpage 

%\begin{center}
 \begin{table}[h!]
% \resizebox{.72\textwidth}{!}{\begin{minipage}{\textwidth}
 \begin{center}
\caption{Final error ${\|A-X\|}_F$ obtained by the algorithms for matrices $A$ of size 10.}
\label{results10}
 \begin{tabular}{|c|c|c|c|c|}
 \hline   &  Type 1 & Type 2 & Type 3 & Type 4    \\
 \hline
 Initial error &   1.50 & 4.16  & 5.25  & 5.24   \\ \hline
 BCD         &   0.74 (170) &   3.37 (170) &   2.36 (179) &   3.68 (162)  \\
 Grad        &   0.87 (114983) &   3.37 (114052) &   2.22 (114701) &   3.66 (114594) \\
FGM   &   0.57 (120641) &   \textbf{3.31} (123055) &   2.07 (118766) &   \textbf{3.46} (116893) \\
 SuccConv    &   \textbf{0.33} (221) &   3.39 (1243) &   1.79 (1272) &   3.48 (1311) \\
 BFGS        &   0.38 (811) &   3.35 (965) &   \textbf{1.61} (1671) &   3.71 (1403) \\ \hline
\end{tabular}
\end{center}
% \end{minipage}}
 \end{table}
 %\end{center}

\newpage 

%\begin{center}
 \begin{table}[h!]
  %\resizebox{.72\textwidth}{!}{\begin{minipage}{\textwidth}
 \begin{center}
\caption{Final error ${\|A-X\|}_F$ obtained by the algorithms for matrices $A$ of size 20.}
\label{results20}
 \begin{tabular}{|c|c|c|c|c|}
 \hline    &  Type 1 & Type 2 & Type 3 & Type 4    \\
 \hline
 Initial error &   2.18 & 6.07  & 9.15  & 10.32   \\  \hline
 BCD         &   1.57 (245) &   4.98 (246) &   3.99 (252) &   6.99 (251)  \\
 Grad        &   1.62 (376522) &   5.02 (373161) &   3.83 (384363) &   6.98 (383315) \\
FGM   &   1.38 (379203) &   \textbf{4.77} (391338) &   3.11 (390444) &   6.57 (388281) \\
 SuccConv    &   1.21 (251) &   5.68 (259) &   2.97 (155) &   \textbf{6.50} (257) \\
 BFGS        &   \textbf{1.18} (1572) &   4.85 (4012) &   \textbf{2.82} (1752) &   6.66 (759) \\ \hline
\end{tabular}
 \end{center}
% \end{minipage}}
 \end{table}
 % \end{center}

\newpage 

%\begin{center}
 \begin{table}[h!]
%   \resizebox{.72\textwidth}{!}{\begin{minipage}{\textwidth}
 \begin{center}
\caption{Final error ${\|A-X\|}_F$ obtained by the algorithms for matrices $A$ of size 50.}
\label{results50}
 \begin{tabular}{|c|c|c|c|c|}
 \hline    &  Type 1 & Type 2 & Type 3 & Type 4    \\
 \hline
 Initial error &   3.50 & 9.77  & 25.36  & 25.90   \\ \hline
 BCD         &   3.03 (133) &   8.56 (133) &   10.79 (139) &   18.19 (140)  \\
 Grad        &   2.90 (121154) &   8.44 (118495) &   10.38 (125581) &   18.28 (125405) \\
 FGM   &   \textbf{2.50} (121385) &   \textbf{8.07} (119355) &   \textbf{8.23} (123023) &   \textbf{17.15} (124506) \\
 SuccConv    &   3.49 (4) &   9.51 (4) &   21.40 (4) &   23.71 (4) \\
 BFGS        &   3.32 (797) &   9.36 (3437) &   19.17 (3116) &   23.93 (435) \\ \hline
\end{tabular}
 \end{center}
 %\end{minipage}}
 \end{table}
 %\end{center}

\newpage 

%\begin{center}
 \begin{table}[h!]
% \resizebox{.72\textwidth}{!}{\begin{minipage}{\textwidth}
 \begin{center}
\caption{Final error ${\|A-X\|}_F$ obtained by the algorithms for matrices $A$ of size 100.}
\label{results100}
 \begin{tabular}{|c|c|c|c|c|}
 \hline   &  Type 1 & Type 2 & Type 3 & Type 4    \\
 \hline
 Initial error &   4.98 & 13.89  & 49.72  & 51.68   \\ \hline
 BCD         &   4.74 (61) &   12.69 (68) &   22.62 (65) &   37.50 (66)  \\
 Grad        &   4.41 (56998) &   12.23 (54735) &   21.34 (59276) &   37.54 (58277) \\
FGM   &   \textbf{3.87} (53768) &   \textbf{11.69} (54603) &   \textbf{16.90} (56733) &   \textbf{35.40} (55890) \\
 SuccConv    &   4.98 (1) &   13.89 (1) &   49.72 (1) &   51.68 (1) \\
 BFGS        &   4.87 (2418) &   13.82 (3316) &   37.61 (1352) &   48.99 (1253) \\ \hline
\end{tabular}
 \end{center}
 % \end{minipage}}
 \end{table}
% \end{center}

\newpage 

We observe the following:
\begin{itemize}
\item In all cases, FGM outperforms BCD and Grad, and is clearly the best to solve our formulation~\eqref{nearportham}.
\item For small matrices ($n \leq 20$), there is no clear winner between FGM, SuccConv and BFGS.
The main reason is that these methods converge to different local minima.
In fact, although they are initialized at the same (stable) starting point (see Section~\ref{init}), they use rather different strategies; FGM uses a gradient-like descent method on problem~\eqref{nearportham}, SuccConv uses succesive convex approximation, and BFGS applies Quasi-Newton to a penalized variant of~\eqref{eq:prob_def}--
so that they can converge to the basin of attraction of different local minima.
As an illustration, Figure~\ref{type1n10} displays the eigenvalues of the nearest stable matrices computed by FGM, SuccConv and BFGS for the matrix of type 1 of size 10.
\begin{figure}[ht]
\begin{center}
\includegraphics[width=0.5\textwidth]{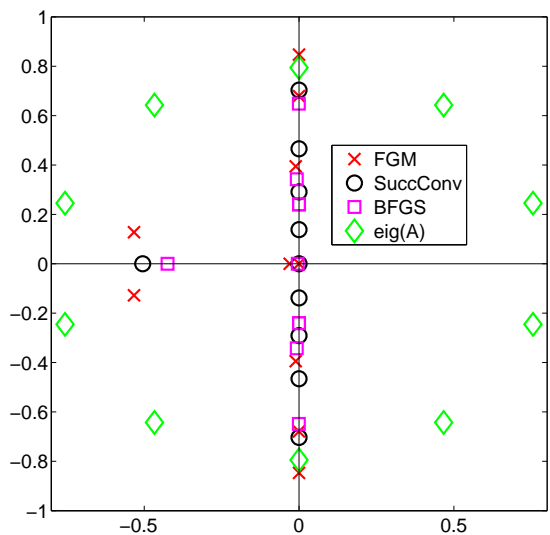}
\caption{Eigenvalues of the solutions computed by FGM, SuccConv and BFGS for the matrix of type 1 of dimension 10.
\label{type1n10}}
\end{center}
\end{figure}
We observe that the solutions are rather different --they correspond to different local minima which explains why the algorithms saturate at different objective function values on Figure~\ref{matn10}.
In particular, FGM has two eigenvalues with large negative real parts,
while the two others algorithms only have one (at different locations).
Note that both solutions of SuccConv and BFGS have two eigenvalues that are located in the right half of the complex plane,
although very close to the imaginary axis (about $10^{-7}$ for SuccConv and $10^{-10}$ for BFGS),
while our algorithm is guaranteed to provide stable matrices (unless 0 is not a semisimple eigenvalue);
see Lemma~\ref{lem:stabJRQ}.

To confirm this behavior, we have run these three algorithms on 100 matrices of type 3 and 4 of size $n=10$.
For matrices of type 3,  FGM provided the best solution for 14 matrices, SuccConv for 46 and BFGS for 40.
FGM (resp.\@ SuccConv and BFGS) had average error 1.95 (resp.\@ 1.83 and 1.86) with standard deviation 0.60 (resp.\@ 0.60 and 0.66).
For matrices of type 4,  FGM provided the best solution for 11 matrices, SuccConv for 35 and BFGS for 54.
FGM (resp.\@ SuccConv and BFGS) had average error 3.30 (resp.\@ 3.28 and 3.26) with standard deviation 0.29 (resp.\@ 0.30 and 0.30).
Therefore, it seems that SuccConv and BFGS tend to identify in general better local minima for these small randomly generated matrices using the initialization from Section~\ref{init},
although it is difficult to explain rigorously.

\item For large matrices ($n \geq 50$), FGM outperforms SuccConv and BFGS providing better approximations much faster.   In particular, SuccConv is extremely slow (performing only 4 iterations within 5 minuts for $n=50$): the reason is that it uses second-order information and each iteration is expensive. BFGS is able to perform more iterations but it is not able to locate good local minima and stabilizes at stationary points with much higher error than BCD, Grad and FGM.

\item Although it is not clear for all the figures, FGM always had the fastest initial convergence, that is, it always decreased the objective function initially the fastest.

\end{itemize}

\section{Conclusion}
We have presented a new way to look at the nearest stable matrix problem by using the concept of dissipative Hamiltonian systems. This resulted in a equivalent optimization problem with a simpler convex feasible set; see the formulation~\eqref{nearportham}.
We have proposed three algorithms namely BCD, Grad and FGM to solve~\eqref{nearportham}. We found that FGM works very well in most situations.
For large matrices ($n \geq 50$), it outperforms all other approaches.
For small matrices ($n \leq 20$), SuccConv and BFGS are sometimes able to identify better local minima starting from the same initial point. \\ We hope this work paved the way to more efficient methods to solve the nearest stable matrix problem based on our new formulation~\eqref{nearportham}. A particular promising way of future research would be to develop more sophisticated techniques to solve~\eqref{nearportham} (e.g., using some second-order information) and to apply globalization schemes to try to find the global minimum.
For example, we observed that using random starting point does not lead to good solutions in general hence coming up with good initialization schemes would also be a direction for further research.

\paragraph{Acknowledgements:}
The authors acknowledge the support of the ERC (starting grant n$^\text{o}$ 679515).
Nicolas Gillis also acknowledges the support of the F.R.S.-FNRS (incentive
grant for scientific research n$^\text{o}$ F.4501.16).

The authors would like to thank Michael Overton and F.-X.\@ Orban de Xivry for providing their code.
They are also grateful to Nicola Guglielmi for helpful discussions.

We also thank the two anonymous reviewers for their insightful comments which helped
improve the paper.

\appendix

\section{Fast gradient method for smooth convex optimization} \label{appA} 

In this section, we briefly describe the fast gradient  method by Nesterov~\cite{nes83}.
We refer the interested reader to the book~\cite{nes04} (in particular, p.87-93). Let us consider the following optimization problem
\[
\min_{x \in \mathcal{X}} f(x),
\]
where $\mathcal{X}$ is a closed convex set, and $f$ is a convex and smooth function with gradient having Lipschitz constant $L$.
If $f$ is quadratic, as it is the case in this paper, $L$ equals the largest eigenvalue of the Hessian $\nabla^2 f$ of $f$. Note that we do not assume here that $f$ is strongly convex.

An optimal method for the problem above is the so-called fast gradient  method (FGM), a.k.a.\@ Nesterov's method.
As opposed to the gradient method that uses one sequence of iterates, FGM uses two sequences which allows to make FGM converge much faster than the standard gradient method (if $f$ is not strongly convex, the difference between the objective function and the optimal value converges at a rate $O(1/k^2)$ for FGM vs.\@ $O(1/k)$ for the gradient method, where $k$ is the iteration number).

Given a point $y$, let us define the gradient mapping
\begin{align*}
x_{\mathcal{X}}(y,L)
& = \argmin_{x \in \mathcal{X}} \left[ f(y) + \langle \nabla f(y), x-y  \rangle + \frac{L}{2} \|x-y\|_2^2 \right] \\
& = \argmin_{x = y+\delta \in \mathcal{X}} \left[ \frac{2}{L} \langle \nabla f(y), \delta  \rangle +  \|\delta\|_2^2 \right] \\
& = \argmin_{x = y+\delta \in \mathcal{X}}  \left\|  \delta +  \frac{1}{L} \nabla f(y) \right\|_2^2  \\
& = \argmin_{x \in \mathcal{X}}  \left\|  x - \left(y - \frac{1}{L} \nabla f(y) \right) \right\|_2^2  \\
& = \mathcal{P}_{\mathcal{X}}  \left(  y - \frac{1}{L} \nabla f(y) \right)  ,
\end{align*}
where $\mathcal{P}_{\mathcal{X}}$ is the projection onto $\mathcal{X}$.
 The standard gradient method will simply apply gradient mapping successively. FGM will use another set of iterates in order to obtain a better convergence rate; see Algorithm~\ref{fastgrad}.

%If the projection can be computed efficiently, then this mapping can be
%Choose the largest steplength so that
%\[
%F (  x_{\mathcal{X}}(y,L)  )
%\leq
%f(y) + \langle \nabla f(y), x_{\mathcal{X}}(y,L) -y  \rangle + \frac{L}{2} ||x_{\mathcal{X}}(y,L) - y||_2^2
%\]

\renewcommand{\thealgorithm}{FGM}
 \algsetup{indent=2em}
\begin{algorithm}[ht!]
\caption{Fast Gradient Method \cite[p.90]{nes04}} \label{fastgrad}
\begin{algorithmic}[1]
\REQUIRE The function $f$ whose gradient has Lipschitz constant $L$,
 an initial guess $x \in \mathcal{X}$, number of iterations $K$ (other stopping criteria can be used).
\ENSURE An approximate solution $x \approx \argmin_{z \in \mathcal{X}} f(z)$.  \medskip
\STATE $\alpha_1 \in (0,1)$; $y = x$ .
\FOR{$k = 1 :$ $K$}
\STATE $\hat x = x$. \hspace{2.55cm} \emph{\% Keep the previous iterate in memory.}
\STATE %$x = \mathcal{P}_{\mathcal{X}}\Big( y - \frac{1}{L} \nabla f(y) \Big)$.
$x = \mathcal{P}_{\mathcal{X}}  \left(  y - \frac{1}{L} \nabla f(y) \right)$ .
\emph{\% $\mathcal{P}_{\mathcal{X}}$ is the projection on $\mathcal{X}$}
\STATE $y = x + \beta_k \left(x - \hat x\right)$, \quad where $\beta_k =  \frac{\alpha_{k} (1-\alpha_{k})}{\alpha_{k}^2 + \alpha_{k+1}}$
with
$\alpha_{k+1} \geq 0$ s.t.
$\alpha_{k+1}^2 = (1-\alpha_{k+1}) \alpha_{k}^2$.
\ENDFOR
\end{algorithmic}
\end{algorithm}

%Try strongly convex ! does not work much better

\small

\bibliographystyle{siam}
\bibliography{bibliostable}

\end{document}